\documentclass[12pt,reqno]{amsart}
\usepackage{amsfonts}
\usepackage{amssymb, latexsym}
\usepackage{eufrak}

\usepackage{amsmath,amsfonts,epsfig}
\usepackage{amsthm}
\usepackage{graphicx}
\usepackage{graphics,color}
\usepackage{amsbsy}
\usepackage{amssymb}

\setbox0=\hbox{$+$}
\newdimen\plusheight
\plusheight=\ht0
\def\+{\;\lower\plusheight\hbox{$+$}\;}

\setbox0=\hbox{$-$}
\newdimen\minusheight
\minusheight=\ht0
\def\-{\;\lower\minusheight\hbox{$-$}\;}

\setbox0=\hbox{$\cdots$}
\newdimen\cdotsheight
\cdotsheight=\plusheight
\def\cds{\lower\cdotsheight\hbox{$\cdots$}}

\renewcommand{\(}{\left\(}
\renewcommand{\)}{\right\)}
\renewcommand{\[}{\left[}

\numberwithin{equation}{section}
 \theoremstyle{plain}
\newtheorem{theorem}{Theorem}[section]
\newtheorem{lemma}[theorem]{Lemma}

\newtheorem{remark}[theorem]{Remark}

\begin{document}
\allowdisplaybreaks
\title[Some results on vanishing coefficients in infinite product expansions] {Some results on vanishing coefficients in infinite product expansions}

\author{Nayandeep Deka Baruah}
\address{Department of Mathematical Sciences, Tezpur University, Napaam-784028, Sonitpur, Assam, INDIA}
\email{nayan@tezu.ernet.in}

\author{Mandeep Kaur}
\address{Department of Mathematical Sciences, Tezpur University, Napaam-784028, Sonitpur, Assam, INDIA}
\email{mandeep@tezu.ernet.in}


\vspace*{0.5in}
\begin{center}
{\bf Some results on vanishing coefficients in infinite product expansions}\\[5mm]
{\footnotesize  NAYANDEEP DEKA BARUAH and MANDEEP KAUR
}\\[3mm]
\end{center}

\vskip 5mm \noindent{\footnotesize{\bf Abstract.}
Recently, M. D. Hirschhorn proved that, if
\begin{align*}
\sum_{n=0}^{\infty}a_nq^n&:=(-q,-q^4;q^5)_\infty(q,q^9;q^{10})_\infty^3\quad \textup{and}\quad\sum_{n=0}^{\infty}b_nq^n:=(-q^2,-q^3;q^5)_\infty(q^3,q^7;q^{10})_\infty^3,
\end{align*}
then $a_{5n+2}=a_{5n+4}=0$ and $b_{5n+1}=b_{5n+4}=0$. Motivated by the work of Hirschhorn, D. Tang proved some comparable results including
the following:

If
\begin{align*}
\sum_{n=0}^{\infty}c_nq^n:=(-q,-q^4;q^5)_\infty^3(q^3,q^7;q^{10})_\infty\quad \textup{and}\quad
\sum_{n=0}^{\infty}d_nq^n:=(-q^2,-q^3;q^5)_\infty^3(q,q^9;q^{10})_\infty,
\end{align*}
then

$c_{5n+3}=c_{5n+4}=0$ and $d_{5n+3}=d_{5n+4}=0$.

In this paper, we prove that
\begin{align*}
a_{5n}&=b_{5n+2},~a_{5n+1}=b_{5n+3},~a_{5n+2}=b_{5n+4},~ a_{5n-1}=b_{5n+1},\\
c_{5n+3}&=d_{5n+3},~ c_{5n+4}=d_{5n+4},~  c_{5n}=d_{5n},~ c_{5n+2}=d_{5n+2},~ \textup{and}~ c_{5n+1}>d_{5n+1}.\end{align*}

 We also record some other comparable results not listed by Tang.

}
\vskip 3mm
\noindent{\footnotesize Key Words: $q$-series expansions; infinite $q$-products; Jacobi's triple product identity, vanishing coefficients.}

\vskip 3mm

\noindent {\footnotesize 2010 Mathematical Reviews Classification
Numbers: Primary 33D15; Secondary 11F33}.

\section{\textbf{Introduction}}\label{secone}
For complex numbers $a$ and $q$, with $|q|<1$, we define
\begin{align*}
(a;q)_\infty&:=\prod_{k=0}^{\infty}(1-aq^k)\\\intertext{and}
(a_1,a_2,\ldots, a_n;q)_\infty&:=(a_1;q)_\infty(a_2;q)_\infty\cdots(a_n;q)_\infty.
\end{align*}

In this paper we prove some new results on vanishing coefficients in the series expansions of certain  infinite $q$-products. In the following few paragraphs we review the work done on this  topic by the previous authors.

In 1978, Richmond and Szekeres \cite{richmond-szekeres} proved that if
\begin{align*}
\sum_{n=0}^{\infty}\alpha_nq^n:=\dfrac{(q^3,q^5;q^8)_\infty}{(q,q^7;q^8)_\infty}\quad \textup{and}\quad
\sum_{n=0}^{\infty}\beta_nq^n:=\dfrac{(q,q^7;q^8)_\infty}{(q^3,q^5;q^8)_\infty},
\end{align*}
then the coefficients $\alpha_{4n+3}$ and $\beta_{4n+2}$ always vanish. They also conjectured that if
\begin{align*}
\sum_{n=0}^{\infty}\gamma_nq^n:=\dfrac{(q^5,q^7;q^{12})_\infty}{(q,q^{11};q^{12})_\infty}\quad \textup{and}\quad
\sum_{n=0}^{\infty}\delta_nq^n:=\dfrac{(q,q^{11};q^{12})_\infty}{(q^5,q^7;q^{12})_\infty},
\end{align*}
then $\gamma_{6n+5}$ and $\delta_{6n+3}$ vanish.

In \cite{andrews-bressoud}, Andrews and Bressoud  proved the following general theorem, which
contains the results of Richmond and Szekeres as special cases.

\begin{theorem}\textup{(Andrews and Bressoud)}\label{theorem1}
If $1\leq r<k$ are relatively prime integers of opposite parity and
$$\dfrac{(q^r,q^{2k-r};q^{2k})_\infty}{(q^{k-r},q^{k+r};q^{2k})_\infty}=:\sum_{n=0}^{\infty}\phi_nq^n,$$
then $\phi_{kn+r(k-r+1)/2}$ is always zero.
\end{theorem}

In \cite{alladi-gordon}, Alladi and Gordon generalized the above theorem as follows:
\begin{theorem}\textup{(Alladi and Gordon)}\label{theorem2}
Let $1<m<k$ and let $(s,km)=1$ with $1\leq s<mk$. Let $r^*=(k-1)s$ and $r\equiv r^*~\textup{mod}~{mk}$ with $1\leq r<mk$.\\
Put $r^\prime =\lceil r^*/{mk}\rceil~\textup{mod}~{k}$ with $1\leq r^\prime<k$. Write
$$\dfrac{(q^r,q^{mk-r};q^{mk})_\infty}{(q^s,q^{mk-s};q^{mk})_\infty}=:\sum_{n=0}^{\infty}\mu_nq^n.$$
Then $\mu_n=0$ for $n\equiv rr^\prime~\textup{mod}~{k}$.
\end{theorem}

They also proved the following companion result to Theorem \ref{theorem2}.

\begin{theorem}\textup{(Alladi and Gordon)} \label{theorem3}
Let $m, k, s, r^*, r$ and $r^\prime$ be defined as in Theorem \ref{theorem2} with $k$ odd. Write
$$\dfrac{(q^r,q^{mk-r};q^{mk})_\infty}{(-q^s,-q^{mk-s};q^{mk})_\infty}=:\sum_{n=0}^{\infty}\mu^\prime_nq^n.$$
Then $\mu^\prime_n=0$ for $n\equiv rr^\prime~\textup{mod}~{k}$.
\end{theorem}

The result of Alladi and Gordon in Theorem \ref{theorem2} does not provide any information about vanishing
coefficients in the cases where $k<m$ or $k = m$. In \cite{mclaughlin}, Mc Laughlin proved the following theorem which covers the cases $k \le m$ as well.

\begin{theorem}\textup{(Mc Laughlin)} \label{theorem4}
Let $k>1$, $m>1$ be positive integers. Let $r=sm+t$, for some integers s and t, where $0\leq s<k$, $1\leq t<m$ and r and k are relatively prime. Let
$$\dfrac{(q^{r-tk},q^{mk-(r-tk)};q^{mk})_\infty}{(q^r,q^{mk-r};q^{mk})_\infty}=:\sum_{n=0}^{\infty}\nu_nq^n;$$
then $\nu_{kn-rs}$ is always zero.
\end{theorem}
He also found the following companion result to Theorem \ref{theorem4}.

\begin{theorem}\textup{(Mc Laughlin)} \label{theorem9}
Let $k>1$, $m>1$ be positive integers with k odd. Let $r=sm+t$, for some integers s and t, where $0\leq s<k$, $1\leq t<m$ and r and k are relatively prime. Let
$$\dfrac{(q^{r-tk},q^{mk-(r-tk)};q^{mk})_\infty}{(-q^r,-q^{mk-r};q^{mk})_\infty}=:\sum_{n=0}^{\infty}\nu^\prime_nq^n;$$
then $\nu^\prime_{kn-rs}$ is always zero.
\end{theorem}

All the proofs of the above theorems use Ramanujan's well-known $_1\psi_1$ summation formula. Very recently, Hirschhorn \cite{hirschhorn} proved the following interesting result by using only the Jacobi triple product identity and elementary $q$-series manipulations.
\begin{theorem}\textup{(Hirschhorn)} \label{hirschhorn}
If
\begin{align*}
\sum_{n=0}^{\infty}a_nq^n&:=(-q,-q^4;q^5)_\infty(q,q^9;q^{10})_\infty^3~ and~\sum_{n=0}^{\infty}b_nq^n:=(-q^2,-q^3;q^5)_\infty(q^3,q^7;q^{10})_\infty^3,
\end{align*}
then \begin{align}\label{hirsch1}a_{5n+2}=a_{5n+4}=0\\\intertext{and} \label{hirsch2}b_{5n+1}=b_{5n+4}=0.\end{align}
\end{theorem}

Note that the forms of the $q$-products in Theorem \ref{hirschhorn} are quite different from those in Theorems \ref{theorem1}--\ref{theorem9}.

Motivated by the work of Hirschhorn \cite{hirschhorn}, Tang \cite{tang} found more results on vanishing coefficients in some other comparable $q$-series expansions. In particular, Tang \cite{tang} proved the following theorem.

\begin{theorem}\textup{(Tang)} \label{tang-thm} If
\begin{align*}
\sum_{n=0}^{\infty}c_nq^n:=(-q,-q^4;q^5)_\infty^3(q^3,q^7;q^{10})_\infty~and~
\sum_{n=0}^{\infty}d_nq^n:=(-q^2,-q^3;q^5)_\infty^3(q,q^9;q^{10})_\infty,
\end{align*}
then
\begin{align}
\label{tang1}c_{5n+3}&=c_{5n+4}=0\\\intertext{and}
\label{tang2}d_{5n+3}&=d_{5n+4}=0.
\end{align}
\end{theorem}

In this paper, we prove the following two theorems. The first theorem clearly implies that instead of proving both \eqref{hirsch1} and \eqref{hirsch2} by Hirschhorn \cite{hirschhorn}, it would have been enough to prove only one of  \eqref{hirsch1} or \eqref{hirsch2}. Similarly, the second theorem implies that instead of proving both \eqref{tang1} and \eqref{tang2} by Tang \cite{tang}, it would have been enough to prove only one of  \eqref{tang1} or \eqref{tang2}.

\begin{theorem}\label{5a}
If
\begin{align*}
\sum_{n=0}^{\infty}a_nq^n&:=(-q,-q^4;q^5)_\infty(q,q^9;q^{10})_\infty^3~ and~\sum_{n=0}^{\infty}b_nq^n:=(-q^2,-q^3;q^5)_\infty(q^3,q^7;q^{10})_\infty^3,
\end{align*}
then \begin{align}\label{ab5n}\sum_{n=0}^{\infty}b_{5n}q^n&-\sum_{n=1}^{\infty}a_{5n-2}q^n=\dfrac{f_1^4}{f_2^4},\\
\label{ab5n1}b_{5n+1}&=a_{5n-1},\\
\label{ab5n2}b_{5n+2}&=a_{5n},\\
\label{ab5n3}b_{5n+3}&=a_{5n+1},\\
\label{ab5n4}b_{5n+4}&=a_{5n+2}.\end{align}
\end{theorem}

\begin{theorem}\label{5b} If
\begin{align*}
\sum_{n=0}^{\infty}c_nq^n:=(-q,-q^4;q^5)_\infty^3(q^3,q^7;q^{10})_\infty~and~
\sum_{n=0}^{\infty}d_nq^n:=(-q^2,-q^3;q^5)_\infty^3(q,q^9;q^{10})_\infty,
\end{align*}
then
\begin{align}
\label{st5n}c_{5n}&=d_{5n},\\
\label{st5n+2}c_{5n+2}&=d_{5n+2},\\
\label{st5n+3}c_{5n+3}&=d_{5n+3},\\\intertext{and}
\label{st5n+4}c_{5n+4}&=d_{5n+4}.
\end{align}
Furthermore,
\begin{align}
\label{st5n+1}\sum_{n=0}^{\infty}c_{5n+1}q^n-\sum_{n=0}^{\infty}d_{5n+1}q^n&=4\dfrac{f_2^4}{f_1^4},
\end{align}
which shows that $c_{5n+1}>d_{5n+1}$.

\end{theorem}

By proceeding in a similar way as in Hirschhorn \cite{hirschhorn}, we also prove the following results.

\begin{theorem}\label{5c}
If
\begin{align*}
(\mp q,\mp q^4;q^5)_\infty(\pm q^4,\pm q^6;q^{10})_\infty^3&=\sum_{n=0}^{\infty}e_nq^n\intertext{and}
(\mp q^2,\mp q^3;q^5)_\infty(\pm q^2,\pm q^8;q^{10})_\infty^3&=\sum_{n=0}^{\infty}f_nq^n,
\end{align*}
where the signs in the products are taken either both upper ones or both lower ones,
then
\begin{align*}
e_{5n+3}=f_{5n+4}=0.
\end{align*}
\end{theorem}

\begin{remark} The results involving the upper ambiguity signs of Theorem \ref{5c} have  already been proved by Tang \cite{tang}. Since our proof works for both the signs, we felt it necessary to keep it here as well.
\end{remark}


\begin{theorem}\label{5e}
If
\begin{align*}
(q,q^4;q^5)_\infty(-q,-q^9;q^{10})_\infty^3&=\sum_{n=0}^{\infty}g_nq^n\intertext{and}
(q^2,q^3;q^5)_\infty(-q^3,-q^7;q^{10})_\infty^3&=\sum_{n=0}^{\infty}h_nq^n,
\end{align*}
then
\begin{align*}
g_{5n+2}=h_{5n+1}=0.
\end{align*}
\end{theorem}

\begin{theorem}\label{5f}
If
\begin{align*}
(q,q^4;q^5)_\infty(q,q^9;q^{10})_\infty^3&=\sum_{n=0}^{\infty}k_nq^n\intertext{and}
(q^2,q^3;q^5)_\infty(q^3,q^7;q^{10})_\infty^3&=\sum_{n=0}^{\infty}\ell_nq^n,
\end{align*}
then
\begin{align*}
k_{5n+4}=\ell_{5n+4}=0.
\end{align*}
\end{theorem}

\begin{theorem}\label{5g}
If
\begin{align*}
(q,q^4;q^5)_\infty^3(-q^3,-q^7;q^{10})_\infty&=\sum_{n=0}^{\infty}s_nq^n\intertext{and}
(q^2,q^3;q^5)_\infty^3(-q,-q^9;q^{10})_\infty&=\sum_{n=0}^{\infty}t_nq^n,
\end{align*}
then
\begin{align*}
s_{5n+3}=t_{5n+4}=0.
\end{align*}
\end{theorem}

\begin{theorem}\label{5h}
If
\begin{align*}
(q,q^4;q^5)_\infty^3(q^3,q^7;q^{10})_\infty&=\sum_{n=0}^{\infty}u_nq^n\intertext{and}
(q^2,q^3;q^5)_\infty^3(q,q^9;q^{10})_\infty&=\sum_{n=0}^{\infty}v_nq^n,
\end{align*}
then
\begin{align*}
u_{5n+4}=v_{5n+3}=0.
\end{align*}
\end{theorem}

We employ simple $q$-series manipulations, Jacobi triple product identity, some preliminary identities for Ramanujan's theta functions, and two known identities for a certain quotient of $q$-products. In Section \ref{sec3}--\ref{sec5}, we prove Theorems \ref{5a}--\ref{5c}, respectively. Since the proofs of Theorems \ref{5e}--\ref{5h} are similar in nature, we omit the proofs.

We now end this section by giving some preliminary results that will be used in our proofs. Let $f(a,b)$ denote  Ramanuajn's theta function defined by
\begin{equation}\label{ab}
f(a,b):=\sum_{k=-\infty}^\infty a^{k(k+1)/2}b^{k(k-1)/2}, \quad \vert ab\vert <1.
\end{equation}
Jacobi's famous triple product identity then takes the form
\begin{align*}
f(a,b)=(-a;ab)_\infty(-b;ab)_\infty(ab;ab)_{\infty}.
\end{align*}

The following preliminary identities easily follow from \cite[p. 46, Entry 30]{bcb3}.
\begin{lemma}\label{lem1} We have \begin{align*}
f(a,ab^2)f(b,a^2b)&=\dfrac{1}{2}f(1,ab)f(a,b),\\
f(a,b)f(-a,-b)&=f(-ab,-ab)f(-a^2,-b^2),\\
f(a,b)&=f(a^3b,ab^3)+af\left(\dfrac{b}{a},a^5b^3\right),\\
f^2(a,b)&=f(a^2,b^2)f(ab,ab)+af\left(\dfrac{b}{a},a^3b\right)f(1,a^2b^2).
\end{align*}
\end{lemma}

The triple product identity and the identities in the above lemma will be used frequently in our proofs, quite often,  without referring.

From \cite[Eqs. (1.19) and (1.20)]{baruah-begum}, we also recall the following two identities which will be used in our next two sections.

\begin{lemma}Let $$R(q)=\dfrac{(q,q^4;q^5)_\infty}{(q^2,q^3;q^5)_\infty}.$$ We have \begin{align*}
\dfrac{1}{R(q)R^2(q^2)}&-q^2R(q)R^2(q^2)=\dfrac{(q^2;q^2)_\infty(q^5;q^5)_\infty^5}{(q;q)_\infty(q^{10};q^{10})_\infty^5}
\\\intertext{and}
\dfrac{R(q^2)}{R^2(q)}-\dfrac{R^2(q)}{R(q^2)}&=4q\dfrac{(q^{10};q^{10})_\infty^5(q;q)_\infty}{(q^5;q^5)_\infty(q^2;q^2)_\infty}.
\end{align*}
\end{lemma}

\section{\textbf{Proof of Theorem \ref{5a}}}\label{sec3}
We have
\begin{align*}
\sum_{n=0}^{\infty}a_nq^n&=(-q,-q^4;q^5)_\infty(q,q^9;q^{10})_\infty^3\\
&=\dfrac{(q^2,q^8;q^{10})_\infty}{(q,q^4;q^5)_\infty}\cdot \dfrac{(q,q^4,q^6,q^9;q^{10})_\infty^3}{(q^4,q^6;q^{10})_\infty^3}\\
&=\dfrac{(q^2,q^8;q^{10})_\infty(q,q^4;q^5)_\infty^2}{(q^4,q^6;q^{10})_\infty^3}\\
&=\dfrac{(q,q^2,q^3,q^4;q^5)_\infty}{(q^2,q^4,q^6,q^8;q^{10})_\infty}\cdot \dfrac{(q,q^4;q^5)_\infty}{(q^2,q^3;q^5)_\infty}\cdot \dfrac{(q^2,q^8;q^{10})_\infty^2}{(q^4,q^6;q^{10})_\infty^2}\\
&=\dfrac{(q;q)_\infty(q^{10};q^{10})_\infty}{(q^2;q^2)_\infty(q^5;q^5)_\infty}\cdot R(q)R^2(q^2)
\end{align*}
and
\begin{align*}
\sum_{n=0}^{\infty}b_nq^n&=(-q^2,-q^3;q^5)_\infty(q^3,q^7;q^{10})_\infty^3\\
&=\dfrac{(q^4,q^6;q^{10})_\infty}{(q^2,q^3;q^5)_\infty}\cdot \dfrac{(q^2,q^3,q^7,q^8;q^{10})_\infty^3}{(q^2,q^8;q^{10})_\infty^3}\\
&=\dfrac{(q^4,q^6;q^{10})_\infty(q^2,q^3;q^5)_\infty^2}{(q^2,q^8;q^{10})_\infty^3}\\
&=\dfrac{(q,q^2,q^3,q^4;q^5)_\infty}{(q^2,q^4,q^6,q^8;q^{10})_\infty}\cdot \dfrac{(q^2,q^3;q^5)_\infty}{(q,q^4;q^5)_\infty}\cdot \dfrac{(q^4,q^6;q^{10})_\infty^2}{(q^2,q^8;q^{10})_\infty^2}\\
&=\dfrac{(q;q)_\infty(q^{10};q^{10})_\infty}{(q^2;q^2)_\infty(q^5;q^5)_\infty}\cdot \dfrac{1}{R(q)R^2(q^2)}.
\end{align*}

Therefore,
\begin{align*}
\sum_{n=0}^{\infty}b_nq^n-\sum_{n=0}^{\infty}a_nq^{n+2}&=\dfrac{(q;q)_\infty(q^{10};q^{10})_\infty}{(q^2;q^2)_\infty(q^5;q^5)_\infty}
\left(\dfrac{1}{R(q)R^2(q^2)}-q^2R(q)R^2(q^2)\right)\\
&=\dfrac{(q;q)_\infty(q^{10};q^{10})_\infty}{(q^2;q^2)_\infty(q^5;q^5)_\infty}\cdot
\dfrac{(q^2;q^2)_\infty(q^5;q^5)_\infty^5}{(q;q)_\infty(q^{10};q^{10})_\infty^5}\\
&=\dfrac{(q^5;q^5)_\infty^4}{(q^{10};q^{10})_\infty^4}.\end{align*}
Equating the coefficients of $q^{5n+r}$, $r=0,1,2,3,4$ from both sides of the above, we readily arrive at
\eqref{ab5n} -- \eqref{ab5n4} to finish the proof.

\section{\textbf{Proof of Theorem \ref{5b}}}\label{sec4}
We have
\begin{align*}
\sum_{n=0}^{\infty}c_nq^n&=(-q,-q^4;q^5)_\infty^3(q^3,q^7;q^{10})_\infty\\
&=\dfrac{(q^2,q^8;q^{10})_\infty^3}{(q,q^4;q^5)_\infty^3}\cdot(q^3,q^7;q^{10})_\infty\\
&=\dfrac{(q^2,q^8;q^{10})_\infty^2(q^2,q^3;q^5)_\infty}{(q,q^4;q^5)_\infty^3}\\
&=\dfrac{(q^2,q^4,q^6,q^8;q^{10})_\infty}{(q,q^2,q^3,q^4)_\infty}\cdot\dfrac{(q^2,q^3;q^5)_\infty^2}{(q,q^4;q^5)_\infty^2}\cdot\dfrac{(q^2,q^8;q^{10})_\infty}{(q^4,q^6;q^{10})_\infty}\\
&=\dfrac{(q^5;q^5)_\infty(q^2;q^2)_\infty}{(q;q)_\infty(q^{10};q^{10})_\infty}\cdot\dfrac{R(q^2)}{R^2(q)}\\\intertext{and}
\sum_{n=0}^{\infty}d_nq^n&=(-q^2,-q^3;q^5)_\infty^3(q,q^9;q^{10})_\infty\\
&=\dfrac{(q^4,q^6;q^{10})_\infty^3}{(q^2,q^3;q^5)_\infty^3}\cdot(q,q^9;q^{10})_\infty\\
&=\dfrac{(q^4,q^6;q^{10})_\infty^2(q,q^4;q^5)_\infty}{(q^2,q^3;q^5)_\infty^3}\\
&=\dfrac{(q^2,q^4,q^6,q^8;q^{10})_\infty}{(q,q^2,q^3,q^4)_\infty}\cdot\dfrac{(q,q^4;q^5)_\infty^2}{(q^2,q^3;q^5)_\infty^2}\cdot
\dfrac{(q^4,q^6;q^{10})_\infty}{(q^2,q^8;q^{10})_\infty}\\
&=\dfrac{(q^5;q^5)_\infty(q^2;q^2)_\infty}{(q;q)_\infty(q^{10};q^{10})_\infty}\cdot\dfrac{R^2(q)}{R(q^2)}.
\end{align*}
Therefore,
\begin{align}\label{5cnew}
\sum_{n=0}^{\infty}c_nq^n-\sum_{n=0}^{\infty}d_nq^n&=\dfrac{(q^5;q^5)_\infty(q^2;q^2)_\infty}{(q;q)_\infty(q^{10};q^{10})_\infty}
\left(\dfrac{R(q^2)}{R^2(q)}-\dfrac{R^2(q)}{R(q^2)}\right)\notag\\
&=4q\dfrac{(q^5;q^5)_\infty(q^2;q^2)_\infty}{(q;q)_\infty(q^{10};q^{10})_\infty}\cdot
\dfrac{(q^{10};q^{10})_\infty^5(q;q)_\infty}{(q^5;q^5)_\infty(q^2;q^2)_\infty}\notag\\
&=4q\dfrac{(q^{10};q^{10})_\infty^4}{(q^5;q^5)_\infty^4}.
\end{align}
Equating the coefficients of $q^{5n+r}$, $r=0,2,3,4$ from both sides of the above, we have
\begin{align*}
c_{5n}&=d_{5n},\\
\sum_{n=0}^{\infty}c_{5n+1}q^n&-\sum_{n=0}^{\infty}d_{5n+1}q^n=4\dfrac{f_2^4}{f_1^4},\\
c_{5n+2}&=d_{5n+2},\\
c_{5n+3}&=d_{5n+3}\\\intertext{and}
c_{5n+4}&=d_{5n+4},
\end{align*}
which are \eqref{st5n} -- \eqref{st5n+4}. Similarly, extracting the terms involving  $q^{5n+1}$ from both sides of \eqref{5cnew}, diving by $q$, and then replacing $q^5$ by $q$, we arrive at \eqref{st5n+1}, to complete the proof.

\section{\textbf{Proof of Theorem \ref{5c}}}\label{sec5}
Throughout this section, we consider the ambiguity signs in the products to be either all upper ones or all lower ones.

We have
\begin{align*}
\sum_{n=0}^{\infty}e_nq^n&=(\mp q,\mp q^4;q^5)_\infty(\pm q^4,\pm q^6;q^{10})_\infty^3\\
&=(\mp q,\mp q^4, \pm q^4, \pm q^4, \pm q^4, \mp q^6, \pm q^6, \pm q^6, \pm q^6, \mp q^9;q^{10})_\infty\\
&=(\mp q,\pm q^4,\pm q^6,\mp q^9;q^{10})_\infty(q^8,q^{12};q^{20})_\infty(\pm q^4;\pm q^6;q^{10})_\infty\\
&=U_1(q)U_2U_3,
\end{align*}
where
$U_1(q)=(\mp q,\pm q^4,\pm q^6,\mp q^9;q^{10})_\infty,~~U_2=(q^8,q^{12};q^{20})_\infty,~~U_3=(\pm q^4;\pm q^6;q^{10})_\infty.$

Now,
\begin{align*}
U_1(-q)&=(\pm q,\pm q^4,\pm q^6,\pm q^9;q^{10})_\infty\\
&=(\pm q, \pm q^4;q^5)_\infty\\
&=\dfrac{(\pm q, \pm q^4,q^5;q^5)_\infty}{(q^5;q^5)_\infty}\\
&=\dfrac{1}{(q^5;q^5)_\infty}\sum_{m=-\infty}^{\infty}(\mp 1)^mq^{(5m^2+3m)/2}\\
&=\dfrac{1}{(q^5;q^5)_\infty}\left(\sum_{m=-\infty}^{\infty}q^{10m^2+3m}\mp q\sum_{m=-\infty}^{\infty}q^{10m^2+7m}\right)\\
&=\dfrac{1}{(q^5;q^{10})_\infty(q^{10};q^{10})_\infty}\left((-q^7,-q^{13},q^{20};q^{20})_\infty\mp q(-q^3,-q^{17},q^{20};q^{20})_\infty\right),
\end{align*}
and hence,
\begin{equation*}
U_1(q)=\dfrac{(q^5;q^5)_\infty (q^{20};q^{20})_\infty }{(q^{10};q^{10})_\infty ^3}\left((q^7,q^{13},q^{20};q^{20})_\infty\pm q(q^3,q^{17},q^{20};q^{20})_\infty\right)
\end{equation*}
Therefore,
\begin{align*}
&U_1(q)U_2\\
&=\dfrac{(q^5;q^5)_\infty}{(q^{10};q^{10})_\infty ^3}(q^8,q^{12},q^{20};q^{20})\big((q^7,q^{13},q^{20};q^{20})_\infty\pm q(q^3,q^{17},q^{20};q^{20})_\infty\big)\\
&=\dfrac{(q^5;q^5)_\infty}{(q^{10};q^{10})_\infty ^3}\sum_{m=-\infty}^{\infty}(-1)^mq^{10m^2+2m}\\
&\quad \times\bigg(\sum_{n=-\infty}^{\infty}(-1)^nq^{10n^2+3n}\pm q\sum_{n=-\infty}^{\infty}(-1)^nq^{10n^2+7n}\bigg)\\
&=\dfrac{(q^5;q^5)_\infty}{(q^{10};q^{10})_\infty ^3}\bigg(\sum_{m,n=-\infty}^{\infty}(-1)^{m+n}q^{10m^2+2m+10n^2+3n}\\
&\quad\pm q\sum_{m,n=-\infty}^{\infty}(-1)^{m+n}q^{10m^2+7m+10n^2+2n}\bigg)\\
&=\dfrac{(q^5;q^5)_\infty}{(q^{10};q^{10})_\infty ^3}\bigg(\bigg(\sum_{r,s=-\infty}^{\infty}q^{10(r+s)^2+2(r+s)+10(r-s)^2+3(r-s)}\\
&\quad-\sum_{r,s=-\infty}^{\infty}q^{10(r+s-1)^2+2(r+s-1)+10(r-s)^2+3(r-s)}\bigg)\\
&\quad \pm q\bigg(\sum_{r,s=-\infty}^{\infty}q^{10(r+s)^2+7(r+s)+10(r-s)^2+2(r-s)}\\&\quad-\sum_{r,s=-\infty}^{\infty}q^{10(r+s-1)^2+7(r+s-1)+10(r-s)^2+2(r-s)}\bigg)\bigg)\\
&=\dfrac{(q^5;q^5)_\infty}{(q^{10};q^{10})_\infty ^3}\bigg(\bigg(\sum_{r,s=-\infty}^{\infty}q^{20r^2+20s^2+5r+s}-q^8\sum_{r,s=-\infty}^{\infty}q^{20r^2+20s^2+15r+21s}\bigg)\\
&\quad \pm q\bigg(\sum_{r,s=-\infty}^{\infty}q^{20r^2+20s^2+9r+5s}-q^3\sum_{r,s=-\infty}^{\infty}q^{20r^2+20s^2+11r+15s}\bigg)\bigg)\\
&=\dfrac{(q^5;q^5)_\infty}{(q^{10};q^{10})_\infty ^3}\bigg((-q^{15},-q^{25},q^{40};q^{40})_\infty\bigg(\sum_{n=-\infty}^{\infty}q^{20n^2+n}\pm q\sum_{n=-\infty}^{\infty}q^{20n^2+9n}\bigg)\\
&\quad \mp q^4(-q^{5},-q^{35},q^{40};q^{40})_\infty\bigg(\sum_{n=-\infty}^{\infty}q^{20n^2+11n}\pm q^4\sum_{n=-\infty}^{\infty}q^{20n^2+21n}\bigg)\bigg).
\end{align*}

We also have
\begin{align*}
U_3&=(\pm q^4,\pm q^6;q^{10})_\infty\\
&=\dfrac{1}{(q^{10};q^{10})_\infty }(\pm q^4,\pm q^6, q^{10};q^{10})_\infty\\
&=\dfrac{1}{(q^{10};q^{10})_\infty }\sum_{m=-\infty}^{\infty}(\mp 1)^mq^{5m^2+m}\\
&=\dfrac{1}{(q^{10};q^{10})_\infty }\bigg(\sum_{m=-\infty}^{\infty}q^{20m^2+2m}\mp q^4\sum_{m=-\infty}^{\infty}q^{20m^2+18m}\bigg).
\end{align*}
It follows that
\begin{align*}
&\sum_{n=0}^{\infty}e_nq^n\\
&=\dfrac{(q^5;q^5)_\infty}{(q^{10};q^{10})_\infty ^4}\bigg(\left(-q^{15},-q^{25},q^{40};q^{40}\right)_\infty \bigg(\sum_{m=-\infty}^{\infty}q^{20m^2+2m}\mp q^4\sum_{m=-\infty}^{\infty}q^{20m^2+18m}\bigg)\\
&\quad \times\bigg(\sum_{n=-\infty}^{\infty}q^{20n^2+n}\pm q\sum_{n=-\infty}^{\infty}q^{20n^2+9n}\bigg)\\
&\quad \mp (-q^{5},-q^{35},q^{40};q^{40})_\infty\bigg(\sum_{m=-\infty}^{\infty}q^{20m^2+2m}\mp q^4\sum_{m=-\infty}^{\infty}q^{20m^2+18m}\bigg)\\
&\quad \times \bigg(q^4\sum_{n=-\infty}^{\infty}q^{20n^2+11n}\pm q^8\sum_{n=-\infty}^{\infty}q^{20n^2+21n}\bigg)\bigg)\\
&=\dfrac{(q^5;q^5)_\infty}{(q^{10};q^{10})_\infty ^4}\big((-q^{15},-q^{25},q^{40};q^{40})_\infty\big(S_1\mp S_2 \pm S_3- S_4\big)\\
&\quad \mp (-q^{5},-q^{35},q^{40};q^{40})_\infty\big(S_5\mp S_6 \pm S_7- S_8\big)\big),
\end{align*}
where
\begin{align*}S_1&=\sum_{m,n=-\infty}^{\infty}q^{20m^2+20n^2+2m+n},\quad S_2=q^4\sum_{m,n=-\infty}^{\infty}q^{20m^2+20n^2+18m+n},\\
S_3&=q\sum_{m,n=-\infty}^{\infty}q^{20m^2+20n^2+2m+9n},\quad
S_4=q^5\sum_{m,n=-\infty}^{\infty}q^{20m^2+20n^2+18m+9n},\\
S_5&=q^4\sum_{m,n=-\infty}^{\infty}q^{20m^2+20n^2+2m+11n},\quad
S_6=q^8\sum_{m,n=-\infty}^{\infty}q^{20m^2+20n^2+18m+11n},\\
S_7&=q^8\sum_{m,n=-\infty}^{\infty}q^{20m^2+20n^2+2m+21n},\quad
S_8=q^{12}\sum_{m,n=-\infty}^{\infty}q^{20m^2+20n^2+18m+21n}.\end{align*}
Proceeding as in \cite{hirschhorn}, it can be shown that the 3-components of the sums $S_1, S_2,\ldots,S_8$ are, respectively,
\begin{align*}&q^{43}\sum_{r,s=-\infty}^{\infty}q^{100r^2+100s^2+125r+40s},\quad q^{23}\sum_{r,s=-\infty}^{\infty}q^{100r^2+100s^2+75r+60s},\\
&q^{23}\sum_{r,s=-\infty}^{\infty}q^{100r^2+100s^2+75r+60s},\quad q^{43}\sum_{r,s=-\infty}^{\infty}q^{100r^2+100s^2+125r+40s},\\
&q^{13}\sum_{r,s=-\infty}^{\infty}q^{100r^2+100s^2+25r+60s},\quad q^{8}\sum_{r,s=-\infty}^{\infty}q^{100r^2+100s^2+25r+40s},\\
&q^{8}\sum_{r,s=-\infty}^{\infty}q^{100r^2+100s^2+25r+40s},\quad q^{13}\sum_{r,s=-\infty}^{\infty}q^{100r^2+100s^2+25r+60s}.\end{align*}
Since these cancel in pairs, we conclude that $e_{5n+3}=0$.

Similarly, we have
\begin{align*}
\sum_{n=0}^{\infty}f_nq^n&=(\mp q^2,\mp q^3;q^5)_\infty(\pm q^2,\pm q^8;q^{10})_\infty^3\\
&=(\mp q^2,\pm q^2, \pm q^2, \pm q^2,\mp q^3 \mp q^7,\mp q^8 \pm q^8, \pm q^8, \pm q^8;q^{10})_\infty\\
&=(\pm q^2,\mp q^3,\mp q^7,\pm q^8;q^{10})_\infty(q^4,q^{16};q^{20})_\infty(\pm q^2;\pm q^8;q^{10})_\infty\\
&=V_1(q)V_2V_3,
\end{align*}
where
$V_1(q)=(\pm q^2,\mp q^3,\mp q^7,\pm q^8;q^{10})_\infty,~~V_2=(q^4,q^{16};q^{20})_\infty,~~V_3=(\pm q^2;\pm q^8;q^{10})_\infty.$

Now
\begin{align*}
V_1(-q)&=(\pm q^2,\pm q^3,\pm q^7,\pm q^8;q^{10})_\infty\\
&=(\pm q^2, \pm q^3;q^5)_\infty\\
&=\dfrac{(\pm q^2, \pm q^3,q^5;q^5)_\infty}{(q^5;q^5)_\infty}\\
&=\dfrac{1}{(q^5;q^5)_\infty}\sum_{m=-\infty}^{\infty}(\mp 1)^mq^{(5m^2+m)/2}\\
&=\dfrac{1}{(q^5;q^5)_\infty}\bigg(\sum_{m=-\infty}^{\infty}q^{10m^2+m}\mp q^2\sum_{m=-\infty}^{\infty}q^{10m^2+9m}\bigg)\\
&=\dfrac{1}{(q^5;q^{10})_\infty(q^{10};q^{10})_\infty}\big((-q^9,-q^{11},q^{20};q^{20})_\infty\mp q^2(-q,-q^{19},q^{20};q^{20})_\infty\big).
\end{align*}
Therefore,
\begin{equation*}
V_1(q)=\dfrac{(q^5;q^5)_\infty (q^{20};q^{20})_\infty }{(q^{10};q^{10})_\infty ^3}\big((q^9,q^{11},q^{20};q^{20})_\infty\mp q^2(q,q^{19},q^{20};q^{20})_\infty\big),
\end{equation*}
and hence,
\begin{align*}
&V_1(q)V_2\\
&=\dfrac{(q^5;q^5)_\infty}{(q^{10};q^{10})_\infty ^3}(q^4,q^{16},q^{20};q^{20})\big((q^9,q^{11},q^{20};q^{20})_\infty\mp q^2(q,q^{19},q^{20};q^{20})_\infty\big)\\
&=\dfrac{(q^5;q^5)_\infty}{(q^{10};q^{10})_\infty ^3}\sum_{m=-\infty}^{\infty}(-1)^mq^{10m^2+6m}\bigg(\sum_{n=-\infty}^{\infty}(-1)^nq^{10n^2+n}\\
&\quad\mp q^2\sum_{n=-\infty}^{\infty}(-1)^nq^{10n^2+9n}\bigg)\\
&=\dfrac{(q^5;q^5)_\infty}{(q^{10};q^{10})_\infty ^3}\bigg(\sum_{m,n=-\infty}^{\infty}(-1)^{m+n}q^{10m^2+6m+10n^2+n}\\
&\quad\mp q^2\sum_{m,n=-\infty}^{\infty}(-1)^{m+n}q^{10m^2+9m+10n^2+6n}\bigg)\\
&=\dfrac{(q^5;q^5)_\infty}{(q^{10};q^{10})_\infty ^3}\bigg(\bigg(\sum_{r,s=-\infty}^{\infty}q^{10(r+s)^2+6(r+s)+10(r-s)^2+(r-s)}\\
&\quad -\sum_{r,s=-\infty}^{\infty}q^{10(r+s-1)^2+6(r+s-1)+10(r-s)^2+(r-s)}\bigg)\\
&\quad \mp q^2\bigg(\sum_{r,s=-\infty}^{\infty}q^{10(r+s)^2+9(r+s)+10(r-s)^2+6(r-s)}\\
&\quad -\sum_{r,s=-\infty}^{\infty}q^{10(r+s-1)^2+9(r+s-1)+10(r-s)^2+6(r-s)}\bigg)\bigg)\\
&=\dfrac{(q^5;q^5)_\infty}{(q^{10};q^{10})_\infty ^3}\bigg(\bigg(\sum_{r,s=-\infty}^{\infty}q^{20r^2+20s^2+7r+5s}-q^4\sum_{r,s=-\infty}^{\infty}q^{20r^2+20s^2+13r+15s}\bigg)\\
&\quad \mp q^2\bigg(\sum_{r,s=-\infty}^{\infty}q^{20r^2+20s^2+15r+3s}-q\sum_{r,s=-\infty}^{\infty}q^{20r^2+20s^2+5r+17s}\bigg)\bigg)\\
&=\dfrac{(q^5;q^5)_\infty}{(q^{10};q^{10})_\infty ^3}\bigg((-q^{15},-q^{25},q^{40};q^{40})_\infty\bigg(\sum_{n=-\infty}^{\infty}q^{20n^2+7n}\pm q^3\sum_{n=-\infty}^{\infty}q^{20n^2+17n}\bigg)\\
&\quad \mp (-q^{5},-q^{35},q^{40};q^{40})_\infty\bigg(q^2\sum_{n=-\infty}^{\infty}q^{20n^2+3n}\pm q^4\sum_{n=-\infty}^{\infty}q^{20n^2+13n}\bigg)\bigg).
\end{align*}

Also,
\begin{align*}
V_3(q)&=(\pm q^2,\pm q^8;q^{10})_\infty\\
&=\dfrac{1}{(q^{10};q^{10})_\infty }(\pm q^2,\pm q^8,q^{10};q^{10})_\infty\\
&=\dfrac{1}{(q^{10};q^{10})_\infty }\sum_{m=-\infty}^{\infty}(\mp 1)^mq^{5m^2+3m}\\
&=\dfrac{1}{(q^{10};q^{10})_\infty }\bigg(\sum_{m=-\infty}^{\infty}q^{20m^2+6m}\mp q^2\sum_{m=-\infty}^{\infty}q^{20m^2+14m}\bigg).
\end{align*}

It follows that
\begin{align*}
&\sum_{n=0}^{\infty}f_nq^n\\
&=\dfrac{(q^5;q^5)_\infty}{(q^{10};q^{10})_\infty ^4}\bigg((-q^{15},-q^{25},q^{40};q^{40})_\infty\bigg(\sum_{m=-\infty}^{\infty}q^{20m^2+6m}\mp q^2\sum_{m=-\infty}^{\infty}q^{20m^2+14m}\bigg)\\
&\quad \times  \bigg(\sum_{n=-\infty}^{\infty}q^{20n^2+7n}\pm q^3\sum_{n=-\infty}^{\infty}q^{20n^2+17n}\bigg)\\
&\quad \mp (-q^{5},-q^{35},q^{40};q^{40})_\infty\bigg(\sum_{m=-\infty}^{\infty}q^{20m^2+6m}\mp q^2\sum_{m=-\infty}^{\infty}q^{20m^2+14m}\bigg)\\
&\quad \times \bigg(q^2\sum_{n=-\infty}^{\infty}q^{20n^2+3n}\pm q^4\sum_{n=-\infty}^{\infty}q^{20n^2+13n}\bigg)\bigg)\\
&=\dfrac{(q^5;q^5)_\infty}{(q^{10};q^{10})_\infty ^4}\bigg((-q^{15},-q^{25},q^{40};q^{40})_\infty\big(T_1\mp T_2 \pm T_3-T_4\big)\\
&\quad \mp (-q^{5},-q^{35},q^{40};q^{40})_\infty\big(T_5\mp T_6 \pm T_7-T_8\big)\big),
\end{align*}
where
\begin{align*}T_1&=\sum_{m,n=-\infty}^{\infty}q^{20m^2+20n^2+6m+7n},\quad
T_2=q^2\sum_{m,n=-\infty}^{\infty}q^{20m^2+20n^2+14m+7n},\\
T_3&=q^3\sum_{m,n=-\infty}^{\infty}q^{20m^2+20n^2+6m+17n},\quad
T_4=q^5\sum_{m,n=-\infty}^{\infty}q^{20m^2+20n^2+14m+17n},\\
T_5&=q^2\sum_{m,n=-\infty}^{\infty}q^{20m^2+20n^2+6m+3n},\quad
T_6=q^4\sum_{m,n=-\infty}^{\infty}q^{20m^2+20n^2+14m+3n},\\
T_7&=q^4\sum_{m,n=-\infty}^{\infty}q^{20m^2+20n^2+6m+13n},\quad
T_8=q^6\sum_{m,n=-\infty}^{\infty}q^{20m^2+20n^2+14m+13n}.\end{align*}

It can be shown that, the 4-components of the sums $T_1, T_2, \ldots ,T_8$ are, respectively,
\begin{align*}&q^{14}\sum_{r,s=-\infty}^{\infty}q^{100r^2+100s^2+20r+75s},\quad
q^{29}\sum_{r,s=-\infty}^{\infty}q^{100r^2+100s^2+80r+75s},\\
&q^{29}\sum_{r,s=-\infty}^{\infty}q^{100r^2+100s^2+80r+75s},\quad
q^{14}\sum_{r,s=-\infty}^{\infty}q^{100r^2+100s^2+20r+75s},\\
&q^{19}\sum_{r,s=-\infty}^{\infty}q^{100r^2+100s^2+80r+25s},\quad
q^{4}\sum_{r,s=-\infty}^{\infty}q^{100r^2+100s^2+20r+25s},\\
&q^{4}\sum_{r,s=-\infty}^{\infty}q^{100r^2+100s^2+20r+25s},\quad
q^{19}\sum_{r,s=-\infty}^{\infty}q^{100r^2+100s^2+80r+25s},\end{align*}
and these cancel in pairs. Therefore, we arrive at $f_{5n+4}=0$ to finish the proof.

\section*{Acknowledgment} The authors would like to thank the anonymous referee for his/her comments. The first author's research was partially supported by Grant no. MTR/2018/000157 of Science \& Engineering Research Board (SERB), DST, Government of India.

\end{document}